# ON DUAL LORENTZIAN HOMOTHETIC EXPONENTIAL MOTIONS WITH ONE PARAMETER

## Mehmet Ali GÜNGÖR, Tülay SOYFİDAN and Muhsin ÇELİK

**Abstract.** In this study, the concept of dual Lorentzian homotetic exponential motions in $\mathbb{D}_1^3$ is discussed and their velocities, accelerations obtained. Also, some geometric results between velocity and acceleration vectors of a point in a spatial motion are obtained. Finally, the theorems related to acceleration and acceleration centres is given.

**AMS 2000 Subject Classification:** 53A17, 53B30, 53B50.

**Keywords:** Dual Lorentz space, homothetic motion, exponential transformation

## 1. Introduction

To investigate the geometry of the motion of a line or a point in the motion of space is important in the study of space kinematics or spatial mechanisms or in physics. The geometry of such a motion of a point or a line has a number of applications in geometric modelling and model-based manufacturing of the mechanical products or in the design of robotic motions. These are specifically used to generate geometric models of shell-type objects and thick surfaces [1, 2, 3].

Dual numbers were introduced in the 19*th* century by Clifford [4] and their first applications to mechanics was generalized by Kotelnikov [5] and Study [6] in their principle of transference. The principle relates spherical and spatial geometry and states as follows; "all laws and formula concerning a spherical configuration (involving intersecting line vectors and real angles) are equally valid when applied to and equivalent spatial configuration of skew vectors, if each real angle $a$, in the original formula is replaced by the corresponding dual angle $a = a + \mathcal{E}a^*$", [7, 8]. In excellence of the principle of transference, spherical geometry can be carried over into the domain of spatial geometry permanently replacing real angles with dual angles and real direction cosines with dual direction cosines, [7, 8].

The dual number has a geometrical meaning which is discussed in detail in [9, 10]. For several decades there were attempts to apply dual numbers to rigid body dynamics. Investigators showed that the momentum of a rigid body can be described as a motor that obeys the motor transformation rule; hence, its derivative with respect to time yields the dual force. However, in those investigations, while going from the velocity motor to the momentum motor, there was always a need to expand the equation to six dimensions and to treat the velocity motor as two separate real vector.

In the Euclidean space $\mathbb{E}^n$ of n-dimensions, W. Clifford and James J. McMahon have given a treatment of a rigid body's motion generated by the most general one parameter affine transformation, [11]. Another treatment was given by H.R. Müller for the same kind of motion, [12]. Subsequently, properties of the planar homothetic motions and three dimensional spherical homothetic motions are given by I. Olcaylar, [13]. The exponential motions were given by A.P. Aydın, [14] and the dual homothetic exponential motions were given by V. Asil, [15].

If we take Minkowski 3-space $\mathbb{R}_1^3$ instead of $\mathbb{E}^3$, the E. Study mapping can be stated as follows: The dual timelike and spacelike unit vectors of dual hyperbolic and Lorentzian unit spheres $\tilde{H}_0^2$ and $\tilde{S}_1^2$ at the dual Lorentzian space $\mathbb{D}_1^3$ are in one-to-one correspondence with the directed timelike and spacelike lines of the space of Lorentzian lines $\mathbb{R}_1^3$, respectively,[16].

In this work, the concept of dual Lorentzian homothetic exponential motions in $\mathbb{D}_1^3$ is studied. Let $\tilde{A}$ $\left(\tilde{A}=-\varepsilon\tilde{A}^T\varepsilon\right)$ and $g \in SO_1(3)$ $\left(g^{-1}=\varepsilon g^T \varepsilon\right)$ represent a constant dual anti-symmetric matrix and a dual orthogonal matrix in dual Lorentz space, respectively. In both representations, $\varepsilon$ is a sign matrix with the diagonal members $\varepsilon_1 = \varepsilon_2 = 1$ and $\varepsilon_3 = -1$. The exponential transformation described as

$$\exp: R \times D_n^n \to \mathrm{GL}(n,D) \subset D_n^n$$

$$(t,\tilde{A}) \to \exp(t\tilde{A}) = g(t) = e^{t\tilde{A}} = \sum_{k=0}^{\infty} \frac{t^k}{k!} \tilde{A}^k$$

is investigated in the view of kinematic under the condition of $\tilde{A}$ $\left(\tilde{A}=-\varepsilon\tilde{A}^T\varepsilon\right)$. It has been shown that the motions corresponding to this mapping are regular in the spacelike and timelike regions, and singular in the null regions. A general formula is obtained for the $n^{th}$ order velocity of this mapping. It has been indicated that, sliding-rolling of the dual pole curves was shown on each other at the motion corresponding this mapping and consequently it was also shown that this type of two dual curves on each other. Then, both the poles of velocity and acceleration belonging to these motions were examined. Finally, the conditions of finding acceleration centres were obtained in the spacelike and timelike regions.

We hope that these results will contribute to the study of space kinematics and physics applications.

**II. Basic Concepts**

If $a$ and $a^*$ are real numbers and $\mathcal{E}^2 = 0$, the combination $\tilde{a} = a + \mathcal{E}a^*$ is called a dual number, where $\mathcal{E}$ is dual unit. The set of all dual numbers forms a commutative ring over the real number field and denoted by $\mathbb{D}$. Then the set

$$D^3 = \left\{ \tilde{a} = (A_1, A_2, A_3) \mid A_i \in D, \ 1 \leq i \leq 3 \right\}$$

is a module over the ring $\mathbb{D}$ which is called a $\mathbb{D}-$module or dual space and denoted by $\mathbb{D}^3$. The elements of $\mathbb{D}^3$ are called dual vectors. Thus a dual vector $\vec{\tilde{a}}$ can be written as

$$\vec{\tilde{a}} = \vec{a} + \mathcal{E}\vec{a^*}$$

where $\vec{a}$ and $\vec{a^*}$ are real vectors in $\mathbb{R}^3$.

The Lorentzian inner product of dual vectors $\vec{\tilde{a}}$ and $\vec{\tilde{b}}$ in $\mathbb{D}^3$ is defined by

$$\langle \vec{\tilde{a}}, \vec{\tilde{b}} \rangle = \langle \vec{a}, \vec{b} \rangle + \mathcal{E}\left( \langle \vec{a}, \vec{b^*} \rangle + \langle \vec{a^*}, \vec{b} \rangle \right)$$

with the Lorentzian inner product $\vec{a}$ and $\vec{b}$ given as

$$\langle \vec{a}, \vec{b} \rangle = a_1 b_1 + a_2 b_2 - a_3 b_3$$

where $\vec{a} = (a_1, a_2, a_3)$ and $\vec{b} = (b_1, b_2, b_3)$. Therefore, $\mathbb{D}^3$ with the Lorentzian inner product $\langle \vec{\tilde{a}}, \vec{\tilde{b}} \rangle$ is called 3-dimensional dual Lorentzian space and denoted by of $\mathbb{D}^3_1$, [16].

A dual vector $\vec{\tilde{a}}$ is said to be timelike if $\vec{a}$ is timelike $(\langle \vec{a}, \vec{a} \rangle < 0)$, spacelike if $\vec{a}$ is spacelike $(\langle \vec{a}, \vec{a} \rangle > 0 \text{ or } \vec{a} = 0)$ and lightlike (or null) if $\vec{a}$ is lightlike $(\langle \vec{a}, \vec{a} \rangle = 0, \ \vec{a} \neq 0)$, where $\langle \ , \ \rangle$ is a Lorentzian inner product with signature $(-,+,+)$. The set of all dual vectors such that $\langle \vec{\tilde{a}}, \vec{\tilde{a}} \rangle = 0$ is called the dual lightlike (or null) cone and is denoted by $\Gamma$. The norm of a dual vector $\vec{\tilde{a}}$ is defined to be

$$\|\vec{\tilde{a}}\| = \|\vec{a}\| + \mathcal{E}\frac{\langle \vec{a}, \vec{a^*} \rangle}{\|\vec{a}\|} = a + \mathcal{E}a^* \ , \ \vec{a} \neq 0.$$

We also consider the time orientation as follows: A dual timelike vector $\vec{\tilde{a}}$ is future–pointing (past-pointing) if and only if $\vec{a}$ is future-pointing (past-pointing).

The dual hyperbolic and dual Lorentzian unit spheres are

$$\tilde{H}^2_0 = \left\{ \vec{\tilde{a}} = \vec{a} + \mathcal{E}\vec{a^*} \in \mathbb{D}^3_1 \mid \|\vec{\tilde{a}}\| = 1, \ \vec{a}, \vec{a^*} \in \mathbb{R}^3_1 \text{ and } \vec{a} \text{ timelike} \right\}$$

and

$$\tilde{S}^2_1 = \left\{ \vec{\tilde{a}} = \vec{a} + \mathcal{E}\vec{a^*} \in D^3_1 \mid \|\vec{\tilde{a}}\| = 1, \ \vec{a}, \vec{a^*} \in R^3_1 \text{ and } \vec{a} \text{ spacelike} \right\},$$

respectively.

The dual Lorentzian cross-product of $\vec{\tilde{a}}$ and $\vec{\tilde{b}}$ is defined as

$$\vec{\tilde{a}} \wedge \vec{\tilde{b}} = \vec{a} \wedge \vec{b} + \mathcal{E}\left(\vec{a} \wedge \vec{b}^* + \vec{a}^* \wedge \vec{b}\right)$$

with the Lorentzian cross-product $\vec{a}$ and $\vec{b}$

$$\vec{a} \wedge \vec{b} = (a_3 b_2 - a_2 b_3, a_3 b_1 - a_1 b_3, a_1 b_2 - a_2 b_1)$$

where $\vec{a} = (a_1, a_2, a_3)$ and $\vec{b} = (b_1, b_2, b_3)$, [16].

### 3. Homothetic Exponential Motions in Dual Lorentz Space

Let $g(t) = e^{t\tilde{A}}$ and $H(t) = \tilde{h}(t) g(t)$ be dual orthogonal matrices in the sense of Lorentz. $\tilde{h}(t)$ is a non-constant dual scalar matrix, $\tilde{A}$ is a constant dual Lorentzian anti-symmetrical matrix (i.e. $\tilde{A}^T = -\varepsilon \tilde{A} \varepsilon$ and $\varepsilon$ is a sign matrix), $t$ a real parameter provided that

$$\begin{bmatrix} \tilde{Y} \\ 1 \end{bmatrix} = \begin{bmatrix} H & \tilde{C} \\ 0 & 1 \end{bmatrix} \begin{bmatrix} \tilde{X} \\ 1 \end{bmatrix}, \qquad (1)$$

which is called a dual Lorentzian homothetic exponential motion in the dual space of 3-dimensions. In equation (1), $\tilde{X}$, $\tilde{Y}$ and $\tilde{C}$ are 3x1 type dual matrices. $g$, $\tilde{h}$ and $\tilde{C}$ are differentiable functions of $C^\infty$ class of the parameter $t$. $\vec{\tilde{X}}$ and $\vec{\tilde{Y}}$ correspond to the position vectors of the same point with respect to the rectangular coordinate frames of the moving space $K_0$ and the fixed space $K$, respectively. At the initial time $t = t_0$ we assume that the coordinate systems in $K_0$ and $K$ are coincident.

To avoid the case of Affine transformation we assume that $\tilde{h} = \tilde{h}(t) \neq \text{constant}$ and to avoid the cases of pure translation and pure rotation we also assume that

$$g' = \tilde{A}g, \; \tilde{C}' \neq 0 \text{ and } H' = \frac{dH}{dt} = \tilde{h}'g + \tilde{h}g' = \left(\tilde{h}' + \tilde{h}\tilde{A}\right)g$$

where $(')$ indicates $d/dt$. On the other hand, since $\tilde{h} = \tilde{h}(t)$ is a scalar matrix, its inverse and transpose are

$$\tilde{h}^{-1} = \frac{1}{\tilde{h}} I \quad \text{and} \quad \tilde{h}^T = \tilde{h},$$

respectively. Since $g$ is a dual Lorentzian orthogonal matrix, the inverse of $H$ is

$$H^{-1} = \tilde{h}^{-1} g^T, \qquad \left(\text{i.e. } g^T = \varepsilon g^{-1} \varepsilon\right).$$

From the equation (1) we can also write

$$\tilde{X} = H^{-1}\tilde{Y} + \tilde{C}_0. \qquad (2)$$

where $-H^{-1}\tilde{C} = \tilde{C}_0$.

Equations (1) and (2) express the coordinate transformations between the fixed and moving space.

**Theorem 1:** In $\mathbb{D}_1^3$ dual Lorentz space, a dual Lorentzian homothetic exponential motion given by equation (1) is not regular in null region.

**Proof:** $H' = \tilde{h}'g + \tilde{h}g' = \left(\tilde{h}' + \tilde{h}\tilde{A}\right)g = \tilde{h}\, g\left(\tilde{A} + \dfrac{\tilde{h}'}{\tilde{h}}I_3\right), \qquad I_3 \in \mathbb{R}_3^3$

where if we define $\tilde{\lambda}(t)$ as $\tilde{\lambda}(t) = -\dfrac{\tilde{h}'(t)}{\tilde{h}(t)}$, then last equation becomes

$$H' = \tilde{h}\, g\left(\tilde{A} - \tilde{\lambda}I_3\right). \qquad (3)$$

From equation (3) we find that

$$\det H' = \det\left(\tilde{h}\,g\right)\det\left(\tilde{A} - \tilde{\lambda}I_3\right) = \mp\tilde{h}^3 \det\left(\tilde{A} - \tilde{\lambda}I_3\right). \qquad (4)$$

As $\det H' = 0$, that is, $H'$ is singular, following equation (4) we get

$$\tilde{h} = 0 \text{ or } \det\left(\tilde{A} - \tilde{\lambda}I_3\right) = 0.$$

Here $\tilde{h} \neq 0$. Otherwise, the motion will be pure translation. Hence $\det\left(\tilde{A} - \tilde{\lambda}I_3\right)$ should be zero. Now we search the solution of equation $\det\left(\tilde{A} - \tilde{\lambda}I_3\right) = 0$. First of all, let us show that

$$\left\langle \tilde{A}\vec{v}, \vec{w} \right\rangle = -\left\langle \vec{v}, \tilde{A}\vec{w} \right\rangle \qquad (5)$$

for $\vec{v}, \vec{w} \in \mathbb{D}_1^3$. Since $\tilde{A}$ is a dual Lorentz anti-symmetric matrix, from equation (5) we see that

$$\left\langle \tilde{A}\vec{v}, \vec{w} \right\rangle = \left(\tilde{A}\vec{v}\right)^T \varepsilon\vec{w} = \vec{v}^T \tilde{A}^T \varepsilon\vec{w} = -\vec{v}^T \varepsilon\tilde{A}\vec{w} = -\left\langle \vec{v}, \tilde{A}\vec{w} \right\rangle.$$

The solution of the equation $\det\left(\tilde{A} - \tilde{\lambda}I_3\right) = 0$ gives the characteristic value $\tilde{\lambda}$ of matrix $\tilde{A}$. We search, for $\vec{x} \in \mathbb{D}_1^3$, the solution $\tilde{A}$ of the following equation

$$\tilde{A}\vec{x} = \tilde{\lambda}\vec{x}. \qquad (6)$$

If we choose $\vec{v} = \vec{w} = \vec{x}$, equations (5) and (6) lead us to

$$\langle \tilde{A}\vec{\tilde{x}}, \vec{\tilde{x}} \rangle = -\langle \vec{\tilde{x}}, \tilde{A}\vec{\tilde{x}} \rangle$$
$$\Rightarrow \langle \tilde{\lambda}\vec{\tilde{x}}, \vec{\tilde{x}} \rangle = -\langle \vec{\tilde{x}}, \tilde{\lambda}\vec{\tilde{x}} \rangle$$
$$\Rightarrow \tilde{\lambda}\langle \vec{\tilde{x}}, \vec{\tilde{x}} \rangle = -\tilde{\lambda}\langle \vec{\tilde{x}}, \vec{\tilde{x}} \rangle$$
$$\Rightarrow 2\tilde{\lambda}\langle \vec{\tilde{x}}, \vec{\tilde{x}} \rangle = 0$$
$$\Rightarrow \tilde{\lambda} = 0 \text{ or } \langle \vec{\tilde{x}}, \vec{\tilde{x}} \rangle = 0.$$

Hence, as $\tilde{\lambda}(t) = -\dfrac{\tilde{h}'(t)}{\tilde{h}(t)}$ and $\tilde{h}(t) \neq \text{constant}$, $\tilde{\lambda}$ cannot be constant and $\langle \vec{\tilde{x}}, \vec{\tilde{x}} \rangle$ must be equal to zero. In this case, if $\langle \vec{\tilde{x}}, \vec{\tilde{x}} \rangle = \langle \vec{x}, \vec{x} \rangle + 2\varepsilon\langle \vec{x}, \vec{x}^* \rangle = 0$, then $\langle \vec{x}, \vec{x} \rangle = 0$. Therefore $\vec{\tilde{x}} \in \mathbb{D}_1^3$ is a null vector and $\tilde{A}\tilde{x} = \tilde{\lambda}\tilde{x}$ has dual solutions with respect to $\tilde{\lambda}$. This at the same time, means that the null vectors should remain invariant under $\tilde{A}$. That is, the null conic remains invariant during the motion. So, except the null vector, $H'$ is always regular. Therefore, we can give the following corollary:

**Corollary 1:** In $\mathbb{D}_1^3$ dual Lorentz space, a dual Lorentzian homothetic exponential motion given by equation (1) is regular in timelike and spacelike regions, and has only one instantaneous rotation centre at all time $t$.

**Theorem 2:** If $H(t) = \tilde{h}(t)g(t)$, $\tilde{h}(t)$ is a dual scalar matrix and $g$ is an dual Lorentzian orthogonal $3 \times 3$ matrix, the $n$ th-order derivatives of $H$ is given by

$$H^{(n)} = \left[\sum_{k=0}^{n} \binom{n}{k} \tilde{h}^{(n-k)} \tilde{A}^k \right] g \ . \qquad (7)$$

**Proof:** The proof of this theorem can be obtained by induction. For $n = 1$,

$$H' = \left(\tilde{h}' + \tilde{h}\tilde{A}\right)g = \left[\sum_{k=0}^{1} \binom{1}{k} \tilde{h}^{(1-k)} \tilde{A}^k \right] g \ .$$

For $n = 2$,

$$H'' = \left(\tilde{h}'' + \tilde{h}'\tilde{A}\right)g + \left(\tilde{h}' + \tilde{h}\tilde{A}\right)\tilde{A}g = \left(\tilde{h}'' + 2\tilde{h}'\tilde{A} + \tilde{h}\tilde{A}^2\right)g$$
$$= \left[\sum_{k=0}^{2} \binom{2}{k} \tilde{h}^{(2-k)} \tilde{A}^k \right] g \ .$$

Hence we have shown that it is true for $(n-1)$. Thus

$$H^{(n-1)} = \left[ \sum_{k=0}^{n-1} \binom{n-1}{k} \tilde{h}^{(n-k-1)} \tilde{A}^k \right] g \ .$$

It can now be shown that it is true for $n$. Thus, for $n$:

$$H^{(n)} = \left[ H^{(n-1)} \right]' = \left\{ \left[ \sum_{k=0}^{n-1} \binom{n-1}{k} \tilde{h}^{(n-k-1)} \tilde{A}^k \right] g \right\}'$$

$$= \left[ \sum_{k=0}^{n-1} \binom{n-1}{k} \tilde{h}^{(n-k)} \tilde{A}^k + \sum_{k=0}^{n-1} \binom{n-1}{k} \tilde{h}^{(n-k-1)} \tilde{A}^{k+1} \right] g \ .$$

By using the Binormal Formula the result can be given as follows:

$$H^{(n)} = \left[ \sum_{k=0}^{n} \binom{n}{k} \tilde{h}^{(n-k)} \tilde{A}^k \right] g \ ,$$

which proves the theorem.

**Theorem 3:** In the $\mathbb{D}_1^3$ dual Lorentz space, the high order velocities of dual Lorentzian homothetic exponential motions are given by

$$\tilde{Y}^{(n)} = \sum_{k=0}^{n} \binom{n}{k} H^{(n-k)} \tilde{X}^{(k)} + \tilde{C}^{(n)}, \qquad \tilde{C}^{(n)} \neq 0 \qquad (8)$$

or

$$\tilde{Y}^{(n)} = \sum_{k=0}^{n} \binom{n}{k} \left[ \sum_{i=0}^{n-k} \binom{n-k}{i} \tilde{h}^{(n-k-i)} \tilde{A}^i \right] g \tilde{X}^{(k)} + \tilde{C}^{(n)} \ . \qquad (9)$$

**Proof:** The proof of this theorem can be obtained by induction. For $n=1$,

$$\tilde{Y}' = H' \tilde{X} + H \tilde{X}' + \tilde{C}' \quad \text{or} \quad \tilde{Y}' = \sum_{k=0}^{1} \binom{1}{k} H^{(1-k)} \tilde{X}^{(k)} + \tilde{C}' \ .$$

For $n=2$,

$$\tilde{Y}'' = H'' \tilde{X} + 2H' \tilde{X}' + H \tilde{X}'' + \tilde{C}'' \quad \text{or} \quad \tilde{Y}'' = \sum_{k=0}^{2} \binom{2}{k} H^{(2-k)} \tilde{X}^{(k)} + \tilde{C}'' \ .$$

Shown that it is true for $n=1$ and $n=2$. We can assume it to be true for $(n-1)$. Now let us show that it is true for $n$. So, for $(n-1)$ let

$$\tilde{Y}^{(n-1)} = \sum_{k=0}^{n-1} \binom{n-1}{k} H^{(n-k-1)} \tilde{X}^{(k)} + \tilde{C}^{(n-1)} \ .$$

Thus, for $n$:

$$\tilde{Y}^{(n)} = \left[\tilde{Y}^{(n-1)}\right]' = \left[\sum_{k=0}^{n-1}\binom{n-1}{k}H^{(n-k-1)}\tilde{X}^{(k)} + \tilde{C}^{(n-1)}\right]'$$

$$= \sum_{k=0}^{n-1}\binom{n-1}{k}H^{(n-k)}\tilde{X}^{(k)} + \sum_{k=0}^{n-1}\binom{n-1}{k}H^{(n-k-1)}\tilde{X}^{(k+1)} + \tilde{C}^{(n)}.$$

By using the Binormal Formula the result can be given as follows:

$$\tilde{Y}^{(n)} = \binom{n}{0}H^{(n)}\tilde{X} + \binom{n}{1}H^{(n-1)}\tilde{X}' + \ldots + \binom{n}{n-1}H'\tilde{X}^{(n-1)} + \binom{n}{n}H\tilde{X}^{(n)} + \tilde{C}^{(n)}$$

or

$$\tilde{Y}^{(n)} = \sum_{k=0}^{n}\binom{n}{k}H^{(n-k)}\tilde{X}^{(k)} + \tilde{C}^{(n)}.$$

From equation (3) we find

$$\tilde{Y}^{(n)} = \sum_{k=0}^{n}\binom{n}{k}\left[\sum_{i=0}^{n-k}\binom{n-k}{i}\tilde{h}^{(n-k-i)}\tilde{A}^i\right]g\,\tilde{X}^{(k)} + \tilde{C}^{(n)}.$$

**Theorem 4:** Let $\tilde{A}$ $\left(\tilde{A} = -\varepsilon\tilde{A}^T\varepsilon\right)$ and $g \in SO_1(3)$ $\left(g^{-1} = \varepsilon g^T\varepsilon\right)$ represent a constant dual anti-symmetric matrix and a dual orthogonal matrix in dual Lorentz space, respectively. The following, items requires existence of each other.

i–)$\tilde{A}g = \tilde{A}$, (10)

ii–)$\tilde{A}^2 = \tilde{A}^3 = \ldots = \tilde{A}^n = 0$. (11)

**Proof:** The proof is clear from the exponential mapping definition given in introduction.

Therefore, according to the Theorem 3 and Theorem 4, if we take the first derivative of $\vec{\tilde{Y}}$, we obtain the following result:

$$\tilde{Y}' = \left(\tilde{h}'g + \tilde{h}\tilde{A}\right)\tilde{X} + \tilde{C}' + \tilde{h}\,g\,\tilde{X}'$$

where $\tilde{V}_a = \tilde{Y}'$ is absolute velocity, $\tilde{V}_f = \left(\tilde{h}'g + \tilde{h}\tilde{A}\right)\tilde{X} + \tilde{C}'$ is the sliding velocity and $\tilde{V}_r = \tilde{h}\,g\,\tilde{X}'$ is the relative velocity of the point $\tilde{X}$ whose position vector is $\vec{\tilde{Y}}$.

Now we look for the point where the sliding velocity of our motion is zero at all time $t$. Such points are called rotational pole or instantaneous rotation centre of the motion. Hence

$$\tilde{V}_f = \left(\tilde{h}'g + \tilde{h}\tilde{A}\right)\tilde{X} + \tilde{C}' = 0. \quad (12)$$

Therefore, according to the Corollary 1 and Theorem 4, equation (12) has only a unique solution. This equation's solution gives us the pole point. Hence, from equation (12) we find

$$\tilde{X} \equiv \tilde{P} = -(H')^{-1} \tilde{C}', \qquad (13)$$

where $\tilde{P}$ is a pole point of moving space. This pole point can be expressed in the fixed system

$$\tilde{Y} \equiv \tilde{Q} = H \tilde{P} + \tilde{C} . \qquad (14)$$

Because points $\tilde{P}$ and $\tilde{Q}$ remain constant at the time $t$ in both systems, these give the equations of fixed and moving pole curves. Differentiating equation (14) with respect to $t$, we obtain $\tilde{Q}' = H' \tilde{P} + H \tilde{P}' + \tilde{C}'$ and using $H' \tilde{P} + \tilde{C}' = 0$, we find

$$\tilde{Q}' = H \tilde{P}' = \tilde{h} g \tilde{P}' . \qquad (15)$$

Equation (15) defines the sliding velocity of point $\tilde{Q}$ at the time $t$. So, we can give the following corollary:

**Corollary 2:** In a dual Lorentzian homothetic exponential motion in $\mathbb{D}_1^3$ dual Lorentz space, tangential vectors of pole curves during motion are coinciding after rotation $g$ and translation $\tilde{h}$.

From equation (15) if we take the norm in the sense of dual Lorentz, we get

$$\|\tilde{Q}'\| = \|H \tilde{P}'\| \Rightarrow \|\tilde{Q}'\| dt = \|\tilde{h} g \tilde{P}'\| dt \Rightarrow \|\tilde{Q}'\| dt = \mp |\tilde{h}|^3 \|\tilde{P}'\| dt, \qquad (16)$$

where $\|\tilde{Q}'\| \neq \vec{0}$ and $\|\tilde{P}'\| \neq \vec{0}$.

If the $d\tilde{s}$ and $d\tilde{s}_1$ are dual arc elements of $(\tilde{P})$ and $(\tilde{Q})$ dual pole curves respectively, then from equation (16) we reach

$$\tilde{s}_1 = \mp \int |\tilde{h}|^3 d\tilde{s} .$$

Therefore we can give the following corollary.

**Corollary 3:** In $\mathbb{D}_1^3$ dual Lorentz space, the pole curves roll with sliding on top of each other in non-null regions through dual Lorentzian homothetic exponential motion. This rolling-sliding motion's coefficients are $\mp \tilde{h}^3$.

## 4. Dual Accelerations and Dual Acceleration Centres

In this section, the accelerations and acceleration centres of a point $\tilde{X}$ at one parameter dual Lorentzian homothetic exponential motion of $K_0$ moving dual Lorentzian space relative to $K$ fixed dual Lorentzian space are investigated.

**Theorem 5:** In $\mathbb{D}_1^3$ dual Lorentz space, there is a relation between the absolute, and the sliding, the relative and Coriolis accelerations of a point $\tilde{X}$ of a moving system under dual Lorentzian homothetic exponential motion as

$$\tilde{\gamma}_a = \tilde{\gamma}_r + \tilde{\gamma}_f + \tilde{\gamma}_c .$$

**Proof:** According to the Theorem 3 and Theorem 4, if we take the second derivative of $\vec{\tilde{Y}}$, we obtain

$$\tilde{Y}'' = \left(\tilde{h}''g + 2\tilde{h}'\tilde{A}\right)\tilde{X} + \tilde{C}'' + 2\left(\tilde{h}'g + \tilde{h}\tilde{A}\right)\tilde{X}' + \tilde{h}\,g\tilde{X}'',$$

where $\tilde{\gamma}_a = \tilde{Y}''$ is absolute acceleration, $\tilde{\gamma}_f = \left(\tilde{h}''g + 2\tilde{h}'\tilde{A}\right)\tilde{X} + \tilde{C}''$ is the sliding acceleration, $\tilde{\gamma}_r = \tilde{h}\,g\,\tilde{X}''$ is the relative acceleration and $\tilde{\gamma}_c = 2\left(\tilde{h}'g + \tilde{h}\tilde{A}\right)\tilde{X}'$ is the Coriolis acceleration of the point $\tilde{X}$ whose position vector is $\vec{\tilde{Y}}$.

Let us assume that $\tilde{\Omega} = 2H\left(\tilde{A} - \tilde{\lambda}I_3\right)H^{-1}$. Taking into consideration of relative velocity and Coriolis acceleration, we get $\tilde{\gamma}_c = \tilde{\Omega}\tilde{V}_r$. From this we give the following corollary.

**Corollary 4:** In $\mathbb{D}_1^3$ dual Lorentz space, Coriolis acceleration vector of a point $\tilde{X}$ in moving space at the time $t$ coincides with the relative velocity vector after the rotation $\tilde{\Omega}$ under dual Lorentzian homothetic exponential motion.

Now, we search the points where sliding accelerations are zero at time $t$.

**Theorem 6:** If a dual Lorentzian homothetic exponential motion given by equation (1) in dual Lorentz space $\mathbb{D}_1^3$ consists of a rotation and a sliding motion, then this dual Lorentzian homothetic exponential motion has always a centre of acceleration when

$$\tilde{h}(t) \neq (c_2 t + c_3) + \varepsilon \left( \frac{c_4 t^2 + c_5 t + c_6}{t + c_0} \right), \quad \tilde{h}(t) \neq (c_7 + c_8 t) e^{-\alpha_1 t} + \varepsilon \left( \frac{c_9 t^3 + c_{10} t^2 + c_{11} t + c_{12}}{t + c_{13}} \right) e^{-\alpha_1 t}$$

and

$$\tilde{h}(t) \neq (c_{14} + c_{15} t) e^{\alpha_1 t} + \varepsilon \left( \frac{c_{16} t^3 + c_{17} t^2 + c_{18} t + c_{19}}{t + c_{20}} \right) e^{\alpha_1 t}.$$

**Proof:** From Theorem (5) we get

$$\tilde{\gamma}_f = (\tilde{h}'' g + 2\tilde{h}' \tilde{A}) \tilde{X} + \tilde{C}'' = 0 \Rightarrow H'' \tilde{X} + \tilde{C}'' = 0 \quad \text{or} \quad H'' \tilde{X} = -\tilde{C}''.$$

From equation (2) we find

$$H'' \tilde{X} = H'' \tilde{C}_0 + 2H' \tilde{C}_0 + H \tilde{C}_0''. \tag{17}$$

For this equation to be solved, $H''$ should be regular. Considering that the matrix $\tilde{A}$ is constant, the derivative of equation (3) gives

$$H'' = H'(\tilde{A} - \tilde{\lambda} I_3) + H(-\tilde{\lambda}' I_3) = H(\tilde{A} - \tilde{\lambda} I_3)(\tilde{A} - \tilde{\lambda} I_3) - H \tilde{\lambda}' I_3 = H\left[ (\tilde{A} - \tilde{\lambda} I_3)^2 - \tilde{\lambda}' I_3 \right]. \tag{18}$$

From the last equation we write

$$\det H'' = \det H \det \left[ (\tilde{A} - \tilde{\lambda} I_3)^2 - \tilde{\lambda}' I_3 \right]$$
$$= \det(\tilde{h} g) \det \left[ (\tilde{A} - \tilde{\lambda} I_3)^2 - \tilde{\lambda}' I_3 \right]$$
$$= \mp \tilde{h}^3 \det \left[ (\tilde{A} - \tilde{\lambda} I_3)^2 - \tilde{\lambda}' I_3 \right].$$

As $\tilde{h} \neq 0$, matrix $H''$ is singular if and if only $\det \left[ (\tilde{A} - \tilde{\lambda} I_3)^2 - \tilde{\lambda}' I_3 \right] = 0$. If we evaluate this determinant we reach

$$\tilde{\mu} \left[ (\tilde{\mu} + \tilde{\alpha}^2)^2 - 4 \tilde{\alpha}^2 \tilde{\lambda}^2 \right] = 0. \tag{19}$$

Here $\vec{\tilde{w}} = \begin{bmatrix} \tilde{w}_1 \\ \tilde{w}_2 \\ \tilde{w}_3 \end{bmatrix}$ is a spacelike or timelike vector and

$$\tilde{\mu} = \tilde{\lambda}^2 - \tilde{\lambda}', \quad \tilde{\alpha}^2 = (\alpha_1 + \varepsilon \alpha_1^*)^2 = \|\vec{\tilde{w}}\|^2. \tag{20}$$

From equation (19) we reach $\tilde{\mu} = 0$ or $(\tilde{\mu} + \tilde{\alpha}^2)^2 - 4 \tilde{\alpha}^2 \tilde{\lambda}^2 = 0$. If $\tilde{\mu} = 0$, from equation (20), that is, $\tilde{\lambda}^2 - \tilde{\lambda}' = 0$ we get

$$\tilde{\lambda} = \frac{-1}{t + c_0} + \varepsilon \frac{c_1}{(t + c_0)^2}, \quad (t \neq -c_0). \tag{21}$$

As $\tilde{\lambda} = -\dfrac{\tilde{h}'}{\tilde{h}}$ then from the last equation we find

$$-\dfrac{\tilde{h}'}{\tilde{h}} = \dfrac{-1}{t+c_0} + \varepsilon \dfrac{c_1}{(t+c_0)^2},$$

$$\tilde{h}(t) = (c_2 t + c_3) + \varepsilon \left( \dfrac{c_4 t^2 + c_5 t + c_6}{t+c_0} \right).$$

where $c_i$ $(0 \leq i \leq 6)$ are arbitrary constants.

If $(\tilde{\mu} + \tilde{\alpha}^2)^2 - 4\tilde{\alpha}^2 \tilde{\lambda}^2 = 0$, then

$$\tilde{\mu} + \tilde{\alpha}^2 = 2\tilde{\alpha}\tilde{\lambda} \text{ or } \tilde{\mu} + \tilde{\alpha}^2 = -2\tilde{\alpha}\tilde{\lambda} .$$

From equation (20), as $\tilde{\mu} = \tilde{\lambda}^2 - \tilde{\lambda}'$, if $\tilde{\mu} + \tilde{\alpha}^2 = 2\tilde{\alpha}\tilde{\lambda}$, we get

$$\tilde{\lambda}^2 - \tilde{\lambda}' + \tilde{\alpha}^2 = 2\tilde{\alpha}\tilde{\lambda}$$

and

$$\tilde{\lambda}' = \tilde{\lambda}^2 - 2\tilde{\alpha}\tilde{\lambda} + \tilde{\alpha}^2 = (\tilde{\lambda} - \tilde{\alpha})^2 .$$

If we substitute $\tilde{\lambda}$ into the last equation, remembering $\tilde{\lambda} = -\dfrac{\tilde{h}'}{\tilde{h}}$, we find

$$\tilde{h}(t) = (c_7 + c_8 t) e^{-\alpha_1 t} + \varepsilon \left( \dfrac{c_9 t^3 + c_{10} t^2 + c_{11} t + c_{12}}{t+c_{13}} \right) e^{-\alpha_1 t}$$

where $c_i$ $(7 \leq i \leq 13)$ are arbitrary constants. Similarly, if $\tilde{\mu} + \tilde{\alpha}^2 = -2\tilde{\alpha}\tilde{\lambda}$ then we find

$$\tilde{h}(t) = (c_{14} + c_{15} t) e^{\alpha_1 t} + \varepsilon \left( \dfrac{c_{16} t^3 + c_{17} t^2 + c_{18} t + c_{19}}{t+c_{20}} \right) e^{\alpha_1 t}$$

where $c_i$ $(14 \leq i \leq 20)$ are arbitrary constants as well. The solutions of equation (19) are

$$\tilde{h}(t) = (c_2 t + c_3) + \varepsilon \left( \dfrac{c_4 t^2 + c_5 t + c_6}{t+c_0} \right) \qquad (\text{for } \tilde{\mu} = 0)$$

and

$$\left.\begin{array}{l} (t) = (c_7 + c_8 t) e^{-\alpha_1 t} + \varepsilon \left( \dfrac{c_9 t^3 + c_{10} t^2 + c_{11} t + c_{12}}{t+c_{13}} \right) e^{-\alpha_1 t} \\[2ex] \tilde{h}(t) = (c_{14} + c_{15} t) e^{\alpha_1 t} + \varepsilon \left( \dfrac{c_{16} t^3 + c_{17} t^2 + c_{18} t + c_{19}}{t+c_{20}} \right) e^{\alpha_1 t} \end{array}\right\} \quad \left( \text{for } (\tilde{\mu} + \tilde{\alpha}^2)^2 - 4\tilde{\alpha}^2 \tilde{\lambda}^2 = 0 \right).$$

This proves the theorem.

*Department of Mathematics, Sakarya University, Sakarya, Turkey*
*E-mail: agungor@sakarya.edu.tr*

*Department of Mathematics, Sakarya University, Sakarya, Turkey*
*E-mail: tsoyfidan@sakarya.edu.tr*

*Department of Mathematics, Sakarya University, Sakarya, Turkey*
*E-mail: mcelik@sakarya.edu.tr*